\documentclass[12pt]{article}
\usepackage{amsmath}
\usepackage{amssymb}
\usepackage{amscd}
\newtheorem{theorem}{Theorem}[section]
\newtheorem{proposition}[theorem]{Proposition}

\newtheorem{lemma}[theorem]{Lemma}
\newtheorem{definition}[theorem]{Definition}

\newtheorem{remark}[theorem]{Remark}

\newenvironment{proof}{\begin{trivlist} \item[]{\bf Proof.}}
{\hfill $\square$\end{trivlist}}
\newcommand{\DSH}{{\rm DSH}}

\renewcommand{\P}{\mathbb{P}}
\newcommand{\C}{\mathbb{C}}

\newcommand{\R}{\mathbb{R}}

\newcommand{\id}{{\rm id}}

\newcommand{\PC}{{\rm PC}}

\newcommand{\dc}{{\rm d^c}}
\newcommand{\ddc}{{\rm dd^c}}
\renewcommand{\d}{{\rm d}}

\newcommand{\loc}{{\rm loc}}
\newcommand{\Lone}{{{\rm L}^1}}
\newcommand{\Ltwo}{{{\rm L}^2}}
\newcommand{\Linfty}{{{\rm L}^\infty}}
\newcommand{\Linftyloc}{{{\rm L}^\infty_\loc}}

\newcommand{\supp}{{\rm supp}}

\title{Dynamics of regular birational maps in $\P^k$}
\author{Tien-Cuong Dinh and Nessim Sibony}
\date{}
\begin{document}
\maketitle
\begin{abstract}
We develop the study of some spaces of currents of bidegree $(p,p)$. 
As an application we
construct the equilibrium measure 
for a large class of birational maps of $\P^k$,  as intersection of Green currents. 
We show that
these currents are extremal and that the corresponding measure is mixing.
\end{abstract}
\section{Introduction}

In \cite{Sibony2}, the second author has introduced a class of
polynomial automorphisms of $\C^k$ -- {\it regular
automorphisms} -- and has constructed for such maps the equilibrium
measures as intersection of invariant positif closed currents --
{\it Green currents} (see also \cite{GuedjSibony, DinhSibony2}).
The measure is proved to be mixing when $k=2$ or $3$.
Regular polynomial automorphisms are Zariski dense in the 
space of  polynomial automorphisms of a given algebraic degree.
In dimension 2, these maps are 
H\'enon type automorphisms (see \cite{BedfordSmillie,
BedfordLyubichSmillie, FornaessSibony1}).

In this paper, we develop the theory of some spaces of currents and we construct
Green currents for a larger class of birational maps of $\P^k$.
We show that the Green currents are
extremal and we obtain a mixing measure as intersection of these
currents. Every small pertubation of regular polynomial
automorphisms belongs to this class. Our method 
can be extended to some rational non-invertible
self-maps of $\P^k$ and to random iteration.

For a H\'enon automorphism $f$ of $\C^2$, it was proved in \cite{FornaessSibony3}
that the Green current $T_+$ is the unique positive closed 
$(1,1)$-current of mass 1 
supported on $K^+:=\{z,\ (f^n(z)) \mbox{ bounded}\}$. In particular, this current
is extremal. The result was extended to regular automorphisms
in \cite{Sibony2} and to weakly regular automorphisms in \cite{GuedjSibony}. Here, we
deal with $(p,p)$-currents, $p>1$. The question is to prove their extremality
which implies the mixing of the equilibrium measure.

The problem was already solved for 
automorphisms of compact K\"ahler manifolds under
the natural assumption that their dynamical degrees are distinct. We proved that 
the Green currents are almost extremal, i.e. 
they belong to finite dimensional extremal faces  
of the cone of positive closed currents. 
We then constructed a mixing measure \cite{DinhSibony4}. 

We use here the same 
{\it method of $\ddc$-resolution} as in 
\cite{DinhSibony1, DinhSibony3, DinhSibony4} to study the Green current of some 
birational maps of $\P^k$. The cohomology space is simpler, but we have to 
extend our calculus to deal with indeterminacy set (see also 
\cite{DinhSibony2, DinhDujardinSibony}).
Most of the paper deals with the extension of the calculus 
to new spaces of currents.
Basically the problem is to
give a meaning to the formula
$$\langle f^*(T),\Phi\rangle= \langle T, f_*(\Phi)\rangle$$
when $f$ has indeterminacy points (see Proposition 3.5). 
We believe that this can be applied in other contexts.

In \cite{Guedj}
Guedj has independently proved, for weakly regular automorphisms of $\C^k$, that
the Green currents of the right degrees are extremal.
\\

We describe now our situation. 
Let $f:\P^k\rightarrow\P^k$ be a birational map of
algebraic degree $d\geq 2$.
Let $I^\pm$ be the indeterminacy set of $f^{\pm 1}$.

\begin{definition} \rm We say that $f$ is {\it regular}
if there exists an integer $s$, $1\leq s\leq k-1$, and
open sets $V^\pm$, $U^\pm$ such that
\begin{enumerate}
\item $\overline V^\pm\cap \overline U^\pm=\emptyset$,
$\overline V^\pm\subset U^\mp$ and $I^\pm\subset V^\pm$.
\item There is a smooth positive closed $(k-s,k-s)$-form $\Theta^+$
supported in $\P^k\setminus\overline V^+$, strictly
positive on $\overline U^+$, and 
a smooth positive closed $(s,s)$-form $\Theta^-$ 
supported in $\P^k\setminus\overline V^-$, strictly
positive on $\overline U^-$.
\item $f$ maps $\P^k\setminus V^+$ into $U^+$;
$f^{-1}$ maps $\P^k\setminus V^-$ into $U^-$.
\end{enumerate}
\end{definition}
Observe that if $\P^k\setminus\overline V^+$
(resp. $\P^k\setminus\overline V^-$) is a union of analytic subsets
of dimension $s$ (resp. $k-s$) of $\P^k$, it carries a form $\Theta^+$
(resp. $\Theta^-$) as above.
If $f$ is regular and $\sigma_1$, $\sigma_2$ 
are automorphisms of $\P^k$ close to the identity,
then $\sigma_1\circ f\circ \sigma_2$ is regular.
When $f$ is a polynomial
automorphism, this definition is equivalent to the definition of
\cite{Sibony2}, i.e. to the fact that $I^+\cap I^-=\emptyset$. 
\\

Consider a regular birational map $f$ of algebraic degree $d\geq 2$.
Let $\delta$ be the algebraic degree of $f^{-1}$.
We show that 
the dynamical degree $d_p$ of $f$ is
equal to $d^p$ for $1\leq p\leq s$, 
the dynamical degree $\delta_q$ of $f^{-1}$ is
equal to $\delta^q$ for $1\leq q\leq k-s$ and $d^s=\delta^{k-s}$
(Proposition 3.2).
We also prove that $f^{\pm 1}$ are
{\it algebraically stable}, i.e. no 
hypersurface is sent under an iterate of $f^{\pm 1}$ to its 
indeterminacy set. Hence, we can construct for $f^{\pm 1}$ Green currents
$T_\pm$ of bidegree $(1,1)$ and of mass $1$. The current $T_+$ (resp. $T_-$)
has H\"older continuous local potentials in $\P^k\setminus
\overline V^+$ (resp. $\P^k\setminus\overline V^-$) and satisfies
the relation $f^*(T_+)=d T_+$ (resp. $f_*(T_-)=\delta T_-$) in $\P^k$
\cite{Sibony2}. 
\\

Let $I_n^\pm$ be the indeterminacy set of $f^{\pm n}$. Define 
$U^+_\infty:=\cup_{n\geq 0} f^{-n}(U^+)\setminus I_n^+$ and 
$U^-_\infty:=\cup_{n\geq 0} f^n(U^-)\setminus I_n^-$.
Our main result is the following theorem.

\begin{theorem} Let $f:\P^k\rightarrow\P^k$ be a regular
birational map as above. Then 
for every $p$, $q$ such that $1\leq p\leq s$ and 
$1\leq q\leq k-s$, the following holds.
\begin{enumerate}
\item[{\rm 1.}] If $T$ is a closed positive $(p,p)$-current 
on $\P^k$ of mass $1$ which belongs to $\PC_p(V^+)$, 
then $d^{-np} f^{n*}(T)$ converge weakly in $U^+_\infty$
to $T^p_+$. If $T$ is a closed positive $(q,q)$-current 
on $\P^k$ of mass $1$ which belongs to $\PC_q(V^-)$,
then $\delta^{-nq} (f^n)_*(T)$ converge weakly in $U^-_\infty$
to $T^q_-$. 
\item[{\rm 2.}] The currents $T^p_+$ and $T^q_-$ are 
extremal in the following sense. For every positive closed
$(p,p)$-current $S$ such that $S\leq T_+^p$ in $\P^k$, 
we have $S=cT_+^p$ in $U^+_\infty$
where $c:=\|S\|$. Analogously for $T_-^q$.
\item[{\rm 3.}] The probability measure $\mu=T^s_+\wedge T^{k-s}_-$ 
is invariant, mixing and 
supported in $U^+\cap U^-$. 
\end{enumerate}
\end{theorem}
The spaces $\PC_p$ will 
be defined in Section 2.
The operator $f^*$ on positive closed currents
will be defined in Section 3. 
We use the method of $\ddc$-resolution (see 
\cite{DinhSibony1, DinhSibony3, DinhSibony4}) in order to prove a
convergence result, stronger than the weak convergence (point 1 of Theorem 1.2).
This will be done in Section 4.
The method gives also a new construction
of Green currents and implies their extremality (point 2
of Theorem 1.2). 
The mixing of $\mu$ is a consequence of point 2 (see 
\cite{Sibony2, GuedjSibony, DinhSibony4} for the proof).

The spaces of currents we use as in \cite{DinhSibony1, DinhSibony3, DinhSibony4} 
are probably of 
interest: they allow to consider intersections of currents 
of bidegree $(p,p)$, $p>1$ (see Remark 2.3).

In \cite{Dinh}, the first named author proved
that $T^s_+$ and $T^{k-s}_-$ are weakly
laminar (see \cite{BedfordLyubichSmillie} for H\'enon maps). 
The H\"older continuity of local potentials of $T_\pm$ on $U^\pm$
implies that the measure $\mu$ is PC. It 
has positive Hausdorff dimension and has no mass on pluripolar
sets (see for example \cite{Sibony2}).   

This article replaces the first version of the same paper of January 2004.

\section{DSH and PC currents}

We will introduce two classes of currents in $\P^k$.
Let $V$ be an open set in $\P^k$. 
The class $\DSH^\bullet(V)$ is {\it the space of test currents}. 
For the bidegree $(0,0)$,
these currents are Differences of q.p.S.H. functions which are pluriharmonic
in a neigbourhood of $\overline V$. 
Recall that an $\Lone$ function $\varphi:\P^k\rightarrow
\R\cup\{-\infty\}$ is {\it q.p.s.h.}
if it is upper semi-continuous and
if $\ddc\varphi\geq -c\omega$, $c>0$, in the sense of currents. 
Here $\omega$ is the standard 
Fubini-Study form on $\P^k$ that we normalize by $\int\omega^k=1$.
A set 
$E\subset \P^k$ is {\it pluripolar} if $E\subset\{\varphi=-\infty\}$ for a
q.p.s.h. function $\varphi$.

The class $\PC_\bullet(V)$ is the space of currents of zero order satisfying 
some regularity property in $\P^k\setminus \overline V$. For example, 
such a positive closed current of bidegree $(1,1)$ has Continuous local Potentials 
in $\P^k\setminus\overline V$ (Proposition 2.2).
\\

Let $\DSH^{k-p}(V)$ denote the space of
real-valued $(k-p,k-p)$-currents $\Phi=\Phi_1-\Phi_2$
on $\P^k$ such that 
\begin{enumerate}
\item $\Phi_i$ are negative, $\Phi_{i|V}$
are $\Linftyloc$ forms on $V$;
\item  $\ddc\Phi_i=\Omega_i^+-\Omega_i^-$ with $\Omega_i^\pm$
positive closed currents
supported in $\P^k\setminus \overline V$.
\end{enumerate}

The mass of a positive or negative current $S$ of bidegree $(k-p,k-p)$ is given 
by the formula $\|S\|:=|\int S\wedge \omega^p|$.
Observe that $\|\Omega_i^+\|=\|\Omega_i^-\|$.
Define
$$\|\Phi\|_\DSH:=
\min\big\{\|\Phi_1\|+\|\Phi_2\|+\|\Omega_1^+\|+\|\Omega_2^+\|,
\Phi_i,\ \Omega_i^\pm \mbox{ as above}\big\}.$$
So, positive closed currents supported in $\P^k\setminus \overline V$
are elements of $\DSH^\bullet(V)$. If $S$ is such a current and $\varphi$ is a 
q.p.s.h. function integrable with respect to the trace measure of $S$, then 
$\varphi S\in\DSH^\bullet(V)$.  

A {\it topology} on $\DSH^\bullet(V)$ is defined as
follows: $\Phi^{(n)}\rightarrow\Phi$ in $\DSH^\bullet(V)$ 
if we can write  $\Phi^{(n)}=\Phi^{(n)}_1-\Phi^{(n)}_2$,
$\ddc\Phi_i^{(n)}=\Omega_i^{(n)+}-\Omega_i^{(n)-}$ as above and 
\begin{enumerate}
\item $\Phi^{(n)}\rightarrow \Phi$ weakly in $\P^k$.
\item $(\|\Phi^{(n)}_i\|+\|\Omega_i^{(n)+}\|)_{n\geq 1}$ is bounded.
\item The $\Phi^{(n)}_i$'s are locally uniformly bounded in $V$.
\item The $\Omega_i^{(n)\pm}$'s are supported in the same compact subset of $\P^k\setminus 
\overline V$.
\end{enumerate}
It is a topology
associated to an inductive limit. 

Observe that smooth forms 
in $\DSH^\bullet(V)$ are dense in 
this space. This can be checked by the standard regularization using automorphisms of $\P^k$.

The following proposition allows to construct currents in
$\DSH^\bullet(V)$ as solutions of $\ddc$-equation 
and shows that they can be used as quasi-potentials of 
positive closed currents (see also \cite{DinhSibony4}). 

\begin{proposition}
Let $\Theta$ be a smooth positive closed $(k-p+1,k-p+1)$-form of mass $1$
supported in a compact $K\subset \P^k\setminus\overline V$.
Let $\Omega$ be a positive closed
$(k-p+1,k-p+1)$-current of mass $m$
supported in $K$. Then,
there exists a \underline{negative} $(k-p,k-p)$-form $\Phi\in
{\cal C}^\infty(\P^k\setminus K)\cap \DSH^{k-p}(V)$
with $\Lone$ coefficients, such that 
$\ddc \Phi = \Omega -m\Theta$. Moreover, $\Phi$ depends linearly and
continuously on $\Omega$. We also have
$\|\Phi\|_{\Linfty(V)}+\|\Phi\|_\DSH\leq c_Km$ where $c_K>0$ is
a constant independent of $\Omega$. The form $\Phi$ is continuous where $\Omega$ is continuous.
\end{proposition}

\begin{proof}
The diagonal $\Delta$ of $\P^k\times\P^k$ is cohomologous to the
positive closed form
$$\alpha(z,w):=\Theta(z)\wedge\omega^{p-1}(w) +\sum_{i\not=k-p+1} \omega^i(z)
\wedge \omega^{k-i}(w).$$
Following \cite[Prop. 6.2.3]{BostGilletSoule}, 
since $\P^k\times\P^k$ is homogeneous, we can find
a negative kernel $G(z,w)$
smooth outside $\Delta$ such that $\ddc G=[\Delta]-\alpha$ and whose coefficients are, in absolute value, 
smaller than
$c|z-w|^{1-2k}$, $c>0$.
Define the negative $\Lone$ form $\Phi$ by
$$\Phi(z):=\int_{w\in\P^k} G(z,w)\wedge \Omega(w).$$
If $\pi_1$ and $\pi_2$ denote the projections of $\P^k\times\P^k$ on its factors,
we have $\Phi=(\pi_1)_*(G\wedge \pi_2^*(\Omega))$ and
$\ddc\Phi=(\pi_1)_*(([\Delta]-\alpha)\wedge \pi_2^*(\Omega))=\Omega-m\Theta$.
The properties of $G$
imply that $\Phi$ is smooth on $\P^k\setminus K$, depends
continuously on $\Omega$ and
$\|\Phi\|_{\Linfty(V)}+\|\Phi\|_\DSH\leq c_Km$. It is clear that $\Phi$ is continuous
where $\Omega$ is continuous.
\end{proof}

Let $\PC_p(V)$ be the space of positive closed $(p,p)$-currents $T$ 
which can be
extended to a linear continuous form on $\DSH^{k-p}(V)$.
The value of this linear
form on $\Phi\in\DSH^{k-p}(V)$ is denoted by $\langle T,\Phi\rangle$. 
Since smooth forms are dense in $\DSH^{k-p}(V)$ the extension is unique.
Of course, if $\ddc\Phi=0$, then $\langle T,\Phi\rangle=\int[T]\wedge [\Phi]$ where $[T]$ and
$[\Phi]$ are classes of $T$ and $\Phi$ in $H^{p,p}(X,\C)$ and $H^{k-p,k-p}(X,\C)$. Indeed, we can 
approach $\Phi$ by $\ddc$-closed forms in $\DSH^{k-p}(V)$ using automorphisms of $\P^k$.

The following proposition justifies our notations which suggest that currents in PC
have some continuity property. Let ${\cal C}_{k-p+1}$
denote the cone of positive closed current $\Omega$ of bidegree $(k-p+1,k-p+1)$
supported in $\P^k\setminus \overline V$. Define a topology on ${\cal C}_{k-p+1}$
as follows: $\Omega_n\rightarrow\Omega$ in ${\cal C}_{k-p+1}$ if 
the $\Omega_n$ are supported in the same compact subset of  $\P^k\setminus \overline V$
and $\Omega_n\rightarrow\Omega$ weakly.

\begin{proposition} Let $T=\alpha+\ddc U$ 
be a positive closed $(p,p)$-current, where
$\alpha$ is a continuous $(p,p)$-form and $U$ is a $(p-1,p-1)$-current on $\P^k$.
\begin{enumerate}
\item[{\rm 1.}] 
If the map $\Omega\mapsto \langle U, \Omega\rangle$, which 
is defined on smooth forms $\Omega\in{\cal C}_{k-p+1}$, can be extended to 
a continuous map on ${\cal C}_{k-p+1}$, then 
$T\in\PC_p(V)$. In particular, if $U$ is a continuous form on 
$\P^k\setminus \overline V$, then $T\in\PC_p(V)$.
\item[{\rm 2.}] 
If $p=1$, then $T\in\PC_1(V)$ if and only if $T$ has Continuous local Potentials 
in $\P^k\setminus\overline V$.
\end{enumerate}
\end{proposition}
\begin{proof} 1. Consider a test current $\Phi\in\DSH^{k-p}(V)$.
Write $\ddc\Phi=\Omega^+-\Omega^-$ where $\Omega^\pm\in {\cal C}_{k-p+1}$.
When $\Phi$ and $\Omega^\pm$ are smooth, we have
$$\langle T, \Phi\rangle =\langle \alpha,\Phi\rangle +\langle U,\ddc \Phi\rangle
=\langle \alpha,\Phi\rangle +\langle U,\Omega^+\rangle - \langle U,\Omega^-\rangle.$$
It is clear that if the map $\Omega\mapsto \langle U, \Omega\rangle$ is 
well defined and continuous on ${\cal C}_{k-p+1}$,
then $\langle T, \Phi\rangle$ 
can be extended to a continuous linear form on $\DSH^{k-p}(V)$. Hence 
$T\in\PC_p(V)$. 

Using Proposition 
2.1, one can prove that the converse is also true. For this, one has only to consider $V$
weakly $(p-1)$-convex (see the definition below) since otherwise the currents in $\DSH^{k-p}(V)$ are $\ddc$-closed.
\\

2. We write $T=\alpha+\ddc U$ with $\alpha$ continuous and $U$ a q.p.s.h. function. 
Let $\Theta$ be a smooth positive $(k,k)$-form of mass 1 supported in 
$\P^k\setminus\overline V$.
Let $a\in \P^k\setminus\overline V$ and $\Phi_a$ be the current 
satisfying $\ddc\Phi_a=\delta_a -\Theta$ given
by Proposition 2.1. When $T\in\PC_1(V)$, using a regularization of $\Phi_a$, we get 
$$\langle T,\Phi_a\rangle = \langle\alpha,\Phi_a\rangle - \langle U,\Theta\rangle 
+U(a).$$
Since $\Phi_a$ and $\langle T,\Phi_a\rangle$ depend continuously on $a$, $U$ is continuous on 
$\P^k\setminus \overline V$.
\par
\end{proof}

\begin{remark} \rm
The notion of PC regularity allows to consider the intersection of currents.
If $T$ belongs to $\PC_p(V)$ and $S$ be a positive closed current supported in 
$\P^k\setminus \overline V$, then the positive closed current $T\wedge S$ is well defined and depends 
continuously on $S$. Indeed, if $\varphi$ is a test real smooth form, $\varphi\wedge S$
belongs to $\DSH^\bullet(V)$. So we can define 
$\langle T\wedge S,\varphi\rangle:=\langle T,\varphi\wedge S\rangle$.
\end{remark}

Assume now that $V$ satisfies some convexity property. We say that $V$ is 
{\it weakly $s$-convex} if there exists a non zero positive closed 
current $\Theta$ of bidegree $(k-s,k-s)$ supported in $\P^k\setminus\overline V$.
By regularization, we can assume that $\Theta$ is smooth. Assume also that $\|\Theta\|=1$.
Observe that every positive closed current of bidegree $(s,s)$ intersects
$\Theta$. Hence, it cannot be supported in $\overline V$.

\begin{proposition} Assume that $V$ is weakly $s$-convex as above.
Let $T\in\PC_p(V)$, $1\leq p\leq s$. There exists $c>0$ such that if $\Phi$ is a negative 
smooth $(s-p,s-p)$-form with $\ddc\Phi\geq -\omega^{s-p+1}$, then
$\langle T,\Phi\wedge\omega^{k-s}\rangle \geq -c(1+\|\Phi\|)$.
In particular,
every q.p.s.h function is integrable
with respect to the trace measure $T\wedge \omega^{k-p}$ and $T$ 
has no mass on pluripolar sets.
\end{proposition}

\begin{proof} 
By scaling, we can assume that $\|\Phi\|\leq 1$.
Hence, $\Phi\wedge\Theta$ belongs to a compact set of $\DSH^{k-p}(V)$.
Since $T$ is in $\PC_p(V)$, there exists $c'>0$
independent of $\Phi$ such that $\langle T,\Phi\wedge\Theta\rangle\geq -c'$. 

On the other hand, if $U$ is a smooth negative $(k-s-1,k-s-1)$-form
such that $\ddc U= \Theta -\omega^{k-s}$, we have
\begin{eqnarray*}
\lefteqn{-\int T \wedge\Phi \wedge \omega^{k-s}+ 
\int T \wedge \Phi \wedge \Theta  = } \\
& = &  
\int T \wedge\Phi \wedge \ddc U = \int T \wedge \ddc\Phi\wedge U
\leq -\int T\wedge \omega^{s-p+1}\wedge U.
\end{eqnarray*} 
We then deduce that $\langle T,\Phi\wedge \omega^{k-s} \rangle\geq -c$
where $c>0$ is independent of $\Phi$.

Now consider a q.p.s.h. function 
$\varphi$ strictly negative on $\P^k$
such that $\ddc\varphi\geq -\omega$.
Let $\varphi_n$ be a sequence of negative smooth functions
decreasing to $\varphi$ such that $\ddc\varphi_n\geq -\omega$.
The first part applied to $\Phi=\varphi_n\omega^{s-p}$ gives
$$\langle T,\varphi_n\omega^{k-p}\rangle \geq -c(1+\|\varphi_n\|_{\Lone})\geq 
-c(1+\|\varphi\|_\Lone).$$
It follows that $\langle T,\varphi\omega^{k-p}\rangle \geq -c(1+\|\varphi\|_\Lone)$.
\end{proof}
The above proposition gives a version of Oka's inequality 
(see \cite{FornaessSibony2}) in the sense that $T$-integrability on the support of 
$\Theta$ implies $T$-integrability.

\begin{proposition} Let $V$ be a weakly $s$-convex open set in $\P^k$ and
$T\in \PC_p(V)$, $1\leq p\leq s-1$. Let
$R$ and $R_i$ be positive closed $(1,1)$-currents. Assume that $R=\omega+\ddc v$
and $R_i=\omega+\ddc v_i$ where $v$ and $v_i$ are continuous
on $\P^k\setminus \overline V$.
Then $R\wedge T$ is well defined and belongs to $\PC_{p+1}(V)$.
In particular,
$R_1\wedge \ldots\wedge R_n$ is well defined and belongs
to $\PC_n(V)$ for $1\leq n\leq s$. If $T_i\rightarrow T$ weakly in 
$\PC_p(V)$ and $v_i\rightarrow v$ locally uniformly
on $\P^k\setminus  \overline V$, 
then $R_i\wedge T_i\rightarrow R\wedge T$ 
weakly in $\PC_{p+1}(V)$.
\end{proposition}

\begin{proof} We can assume that 
$v$ is negative.
Proposition 2.4 permits to
define $R\wedge T:=\omega\wedge T +\ddc (vT)$ 
(even without assuming that $v$ is continuous).
It is easy to check by approximation that $R\wedge T$ is positive.
If $\Phi\in\DSH^{k-p-1}(V)$ is a smooth form, we have
$$\langle R\wedge T, \Phi\rangle := \langle T,
\omega\wedge\Phi \rangle + \langle T, v \ddc \Phi\rangle.$$  
When $\Phi\in \DSH^{k-p-1}(V)$ is not smooth,
the right hand side is well defined and depends
continuously on $\Phi$ (see Remark 2.3 for the definition of the measure $T\wedge
\ddc\Phi$). 
Hence, we can extend $R\wedge T$ to a linear continuous form on 
$\DSH^{k-p-1}(V)$.
It follows that $R\wedge
T\in\PC_{p+1}(V)$.
For the second part of Proposition 2.5, it follows from
Proposition 2.2 that $R_1\in\PC_1(V)$.
We then use an induction on $n$.

To prove the convergence result, we use the above formula:
$$\langle R_i\wedge T_i, \Phi\rangle := \langle T_i,
\omega\wedge\Phi \rangle + \langle T_i, v_i \ddc \Phi\rangle.$$
The convergence of the first term is clear for $\Phi\in\DSH^{k-p-1}(V)$. For
the second term, observe that $T_i\wedge\ddc\Phi$ 
are measures with bounded mass supported in the same compact 
subset of $\P^k\setminus\overline V$.
The convergence follows.
\end{proof}

Let $V$ be as in Proposition 2.5 and let $A$ be 
a compact analytic subset of $\P^k\setminus \overline V$. 
Define ${\cal C}$ the cone
of negative $\Lone$ forms 
$\Phi\in {\cal C}^0(\P^k\setminus A)\cap \DSH^{k-p}(V)$
such that $\ddc \Phi=\Omega^+-\Omega^-$
with $\Omega^\pm$ positive closed supported in $\P^k\setminus \overline V$,
continuous on $\P^k\setminus A$
and having no mass on $A$. Here, ${\cal C}^0(\P^k\setminus A)$ denotes the space of continuous 
forms on $\P^k\setminus A$. 
We will use the following lemma in Section 4.
\begin{lemma} Let $R_i$ as in Proposition 2.5. 
Let $S$ be a positive closed $(p,p)$-current, $1\leq p\leq s$, 
such that $S\leq R_1\wedge\ldots\wedge R_p$. 
Then
$S$ can be extended to a continuous linear form on ${\cal C}$ by
$$\langle S,\Phi\rangle := \langle S,\Phi\rangle_{\P^k\setminus A}:=\int_{\P^k\setminus A} S\wedge \Phi.$$
The continuity is with respect to 
the topology of ${\cal C}^0(\P^k\setminus A)\cap \DSH^{k-p}(V)$.  
\end{lemma}

\begin{proof} Define $T_i:=R_1\wedge\ldots\wedge R_i$.
Let $\Phi_n\in\DSH^{k-p}(V)$ be smooth negative forms on
$\P^k$ such that $\Phi_n\rightarrow \Phi$ in
${\cal C}^0(\P^k\setminus A)$ and in $\DSH^{k-p}(V)$.
We show that $\lim \langle S,\Phi_n\rangle = \langle
S,\Phi\rangle_{\P^k\setminus A}$. This will prove the Lemma.  
We have by Fatou's lemma:
$$\limsup \langle S,\Phi_n\rangle  \leq \langle
S,\Phi\rangle_{\P^k\setminus A}$$
and 
$$\limsup \langle T_p-S,\Phi_n\rangle \leq \langle
T_p-S,\Phi\rangle_{\P^k\setminus A}.$$
Since, by Proposition 2.5,  $\langle T_p,\Phi_n\rangle \rightarrow \langle T_p,\Phi\rangle$, 
we only need to
prove that $\langle T_p,\Phi\rangle = \langle
T_p,\Phi\rangle_{\P^k\setminus A}$. 
 
Let $u$ be a negative q.p.s.h. function such that $\ddc u\geq
-\omega$, $u=-\infty$ on $A$ and $u$ is smooth on $\P^k\setminus A$. 
Let $\chi$ be a smooth convex
increasing function on $\R^-\cup\{-\infty\}$ such that $\chi(0)=1$, $\|\chi\|_{{\cal
C}^2}\leq 4$ and $\chi=0$
on $[-\infty,-1]$.
Define $u_n:=\chi(u/n)$. These functions are smooth, equal to $0$
in neigbourhoods of $A$.
We also have $\ddc u_n\geq -4n^{-1}\omega$ and
$u_n\rightarrow 1$ uniformly on compact sets of $\P^k\setminus A$. 
It is sufficient to show that 
$\lim \langle T_p,u_n\Phi\rangle = \langle T_p,\Phi\rangle$.

Let $\ddc\Phi=\Omega^+-\Omega^-$ and define $\Omega:=\Omega^+-\Omega^-$. 
We have 
$$\langle T_p,\Phi\rangle = \langle v_p T_{p-1}, \Omega\rangle + 
\langle T_{p-1}, \omega\wedge \Phi\rangle.$$
This is true for smooth forms and hence for $\Phi$ 
by approximation. On the other hand, we have
$$\langle T_p,u_n\Phi\rangle = \langle \ddc v_p \wedge T_{p-1}, u_n\Phi\rangle + 
\langle T_{p-1}, \omega\wedge u_n\Phi\rangle.$$
Using an induction on $p$, we only need to prove that
$$\lim \langle \ddc v_p\wedge T_{p-1}, u_n\Phi\rangle = 
\langle v_p T_{p-1}, \Omega\rangle.$$ 

Let $\epsilon>0$, $U\Subset \P^k\setminus \overline V$ be a neighbourhood of $A$, $M$
a constant such that $M\geq -\inf_U v_p$, and $v_p^M:=\max(v_p,-M)$. Since $v_p$ is continuous on 
$\P^k\setminus \overline V$, $\ddc v_p^M\wedge T_{p-1}\wedge\Phi\rightarrow \ddc v_p\wedge T_{p-1}\wedge \Phi$
on $\P^k\setminus A$. The measures $\ddc v_p^M\wedge T_{p-1}\wedge \Phi$ and $\ddc v_p\wedge T_{p-1}\wedge \Phi$
are equal in $U$. Since $u_n\rightarrow 1$ locally uniformly on $\P^k\setminus A$,
there exists $n_0$ such that if $n\geq n_0$ we have
$$|\langle \ddc v_p\wedge T_{p-1}, u_n\Phi\rangle - \langle \ddc v^M_p\wedge T_{p-1}, u_n\Phi\rangle |
\leq\epsilon.$$
Hence, if we replace $v_p$ by $v_p^M+M$, we can assume that $v_p$ is positive. In particular, 
$v_p^2$ is q.p.s.h. Hence, $\ddc(v_p^2 T)$ is a difference of positive closed currents. It follows that
$\d v_p$ and $\dc v_p$ belong to $\Ltwo(T_{p-1})$.

We have
\begin{eqnarray*}
\langle \ddc v_p\wedge T_{p-1}, u_n\Phi\rangle
& = & \langle u_nv_p T_{p-1}, \Omega\rangle
-\langle \d u_n\wedge\dc v_p\wedge T_{p-1}, \Phi\rangle + \\
& & 
+ \langle \dc u_n\wedge\d v_p\wedge T_{p-1}, \Phi\rangle
- \langle \ddc u_n\wedge v_p T_{p-1}, \Phi\rangle.
\end{eqnarray*}
By induction hypothesis, the measure $T_{p-1}\wedge\Omega$ has no mass
on $A$ (see also Remark 2.3). Hence, the first term tends to 
$\langle v_pT_{p-1},\Omega\rangle$. 
We show that the other terms
tend to 0. 

Since $\pm\ddc u_n\leq \ddc u_n+8n^{-1}\omega$ and 
$\ddc u_n+8n^{-1}\omega\geq 0$, we have:
\begin{eqnarray*}
|\langle \ddc u_n\wedge v_p T_{p-1}, \Phi\rangle| & \lesssim &
- \langle \ddc u_n \wedge T_{p-1}+ 8 n^{-1}\omega\wedge T_{p-1}, 
\Phi \rangle \\
& \lesssim & 
-\langle  T_{p-1}, u_n\ddc\Phi\rangle - 8n^{-1}\langle T_{p-1}, 
\omega\wedge \Phi \rangle.
\end{eqnarray*}
It follows that  $\langle \ddc u_n\wedge v_p T_{p-1}, \Phi\rangle$
tends to $0$. Indeed, since $u_n\ddc\Phi\rightarrow \ddc\Phi$ in $\DSH^{k-p+1}(V)$
and $T_{p-1}\in \PC_{p-1}(V)$, we have $\langle T_{p-1}, u_n\ddc\Phi\rangle\rightarrow \langle T_{p-1},\ddc\Phi
\rangle=0$. 

For the other terms it is sufficient to use 
the Cauchy-Schwarz inequality and the property that $\d u_n\wedge \dc u_n$
can be dominated by
$\ddc u_n^2+ 100 n^{-1}\omega$. The functions $u_n^2$ satisfy analogous 
inequalities as the $u_n$ do.
\end{proof}

\section{Regular birational maps}

Let $f:\P^k\rightarrow\P^k$ be a dominating rational map of
algebraic degree $d\geq 2$.
In homogeneous coordinates $[z_0:\cdots:z_k]$, we have
$f=[P_0:\cdots:P_k]$ where $P_i$ are homogeneous polynomials of
degree $d$ without common divisor. 
Let $\Gamma$ be the graph of $f$ in $\P^k\times \P^k$,
$\pi_i$ the canonical projections of $\P^k\times\P^k$ onto its factors.
If $A$ is a subset of $\P^k$, define
$f(A):=\pi_2\big(\pi_1^{-1}(A)\cap\Gamma\big)$ and
$f^{-1}(A):= \pi_1\big(\pi_2^{-1}(A)\cap\Gamma\big)$.
The operators $f_*:=(\pi_2)_*(\pi_{1|\Gamma})^*$ and
$f^*:=(\pi_1)_*(\pi_{2|\Gamma})^*$ are well defined and continuous
on $\Linfty$ forms (forms with $\Linfty$ coefficients) with value in spaces of $\Lone$ forms
(forms with $\Lone$ coefficients). 

We define the {\it dynamical degree of order $p$} of
$f$ by
\begin{eqnarray}
d_p & := & \lim_{n\rightarrow\infty} \|f^{n*}(\omega^p)\|^{1/n}
=\lim_{n\rightarrow\infty} \left(\int_{\P^k} f^{n*}(\omega^p)\wedge
\omega^{k-p}\right)^{1/n} \nonumber \\
& = & 
\lim_{n\rightarrow\infty} \|(f^n)_*(\omega^{k-p})\|^{1/n}=
\lim_{n\rightarrow\infty}
\left(\int_{\P^k} (f^n)_*(\omega^{k-p})\wedge
\omega^p\right)^{1/n}
\end{eqnarray} 
These limits always exist \cite{DinhSibony5}.
It is easy to see that $d_p\leq d_1^p$.
The last degree $d_k$ is {\it the topological degree} of $f$.
It is equal to $\#f^{-1}(z)$ for $z$ generic.
\\
 
Consider now, a
birational map $f$, i.e. a map with topological degree 1.
The set $I^+$ (resp. $I^-$) of points $z\in
\P^k$ such that $f(z)$ (resp. $f^{-1}(z)$) is infinite is the
{\it indeterminacy set} of $f$ (resp. $f^{-1}$). 
Hence $f\circ f^{-1}=f^{-1}\circ f=\id$ out of an analytic set.
Let $\delta$ denote the algebraic
degree and $\delta_q$ the
dynamical degree of order $q$ associated to $f^{-1}$. 

\begin{definition} \rm We say that $f$ is {\it $s$-regular},
$1\leq s\leq k-1$, if there exist two 
open sets $V$, $U$ such that
\begin{enumerate}
\item $\overline V\cap \overline U=\emptyset$,
$I^+\subset V$ and $I^-\subset U$.
\item There is a smooth positive closed $(k-s,k-s)$-form $\Theta$
supported in $\P^k\setminus\overline V$ and strictly
positive on $\overline U$. We will assume that
$\|\Theta\|=1$.
\item $f$ maps $\P^k\setminus V$ into $U$.
\end{enumerate}
\end{definition}
Observe that $V$ is weakly $s$-convex. 
If $H$ is a hypersurface of $\P^k$, then $H\not\subset
\overline V$. It follows that $H$ cannot be sent by an iterate
of $f$ to $I^+$. Hence, $f$ is algebraically stable, i.e. $\deg (f^n)=d^n$ 
\cite{Sibony2}.

\begin{proposition} Let $f:\P^k\rightarrow\P^k$ be an $s$-regular
birational map as in Definition 3.1.
Let $I_n^\pm$ be the indeterminacy set of $f^{\pm n}$.
Then $I_n^+\subset V$, $I_n^{-}\subset U$, $\dim I_n^+\leq k-s-1$
and $d_p=d^p$ for $1\leq p\leq s$. 
We have $(f^n)^*=(f^*)^n$ on $H^{p,p}(X,\C)$ for $1\leq p\leq s$. 
If $f$ is regular
as in Definition 1.1, 
then $\dim I_n^-\leq s-1$, 
$\delta_q=\delta^q$ for $1\leq q\leq k-s$ and $d^s=\delta^{k-s}$.
\end{proposition}

\begin{proof} Since $f^n$ is holomorphic on a neighbourhood of
$\P^k\setminus V$, we have $I^+_n\subset V$.
Since $f^{-1}:\P^k\setminus U\rightarrow V$ 
is holomorphic, we have $I_n^-\subset U$.
If $\dim I_n^+\geq k-s$, then
the current
of integration on $I_n^+$ intersects
$\Theta$ which is cohomologous to $\omega^{k-s}$
(recall $\dim H^{p,p}(\P^k,\C)=1$). This is impossible since
$I^+_n\subset V$ and $\supp(\Theta)\cap V=\emptyset$.

Since $f$ is algebraically stable, 
$f^{n*}(\omega)$ is a positive closed $(1,1)$-current of mass
$d^n$ and smooth on $\P^k\setminus I^+_n$. 
We have seen that $\dim I^+_n\leq k-s-1$. The intersection theory 
\cite{Demailly, FornaessSibony2} implies that
$f^{n*}(\omega)\wedge \ldots \wedge f^{n*}(\omega)$ ($p$ times, $p\leq s$)
is well defined and does not charge algebraic sets. It's mass is equal to $d^{np}$.
We deduce from (1) that $d_p=d^p$ and $(f^n)^*=(f^*)^n$ on $H^{p,p}(\P^k,\C)$.

When $f$ is regular, we prove in
the same way that $\dim I_n^-\leq s-1$ and
$\delta_q=\delta^q$. We obtain from (1) that $d_s=\delta_{k-s}$. 
It follows that $d^s=\delta^{k-s}$.
\end{proof}
\begin{remark} \rm
The identity $(f^n)^*=(f^*)^n$ 
on $H^{p,p}(X,\C)$ corresponds to 
an algebraic stability of higher order. 
The notion can be introduced for meromorphic maps on a compact K\"ahler manifold.
Proposition 3.2 is valid in a more general 
case.
\end{remark}

Let $T$ be a positive closed $(p,p)$-current on $\P^k$.
The restriction $f_0$ of $f$ to $\P^k\setminus f^{-1}(I^-)\cup I^+$ is an injective holomorphic map.
We can define $f_0^*(T)$
on $\P^k\setminus f^{-1}(I^-)\cup I^+$. By approximation, one can check that this is a positive
closed current of finite mass (see also \cite{DinhSibony5}). 
Let $f^\star(T)$ denote the trivial extension of $f_0^*(T)$ on $\P^k$. By a
theorem of Skoda \cite{Skoda},
$f^\star(T)$ is positive and closed. 

If $T_n\rightarrow
T$, we have $f^\star (T_n)\rightarrow f^\star (T)$ on
$\P^k\setminus f^{-1}(I^-)\cup I^+$. Moreover, $f^\star(T)$ is smaller than
every limit value $\tau$ of the sequence $f^\star(T_n)$. 
More precisely, the current $\tau-f^\star(T)$ is positive closed and supported 
in $f^{-1}(I^-)\cup I^+$. 

Assume now that $1\leq p\leq s$. Proposition 3.2 implies that $\|f^*(T)\|=d^p\|T\|$ for $T$ smooth.
Using a regularization of $T$, we deduce from the above properties 
that $\|f^\star(T)\| \leq
d^p\|T\|$. 
When $\|f^\star(T)\| = d^p\|T\|$, we define
$f^*(T):=f^\star(T)$. 
We define similarly $f_\star$ and $f_*$ on positive
closed currents.

\begin{lemma}
The operator $f^*$ is continuous: 
if $f^*(T_n)$ and $f^*(T)$ are well 
defined in the above sense 
and if $T_n\rightarrow T$ then $f^*(T_n)\rightarrow f^*(T)$.
If $f^*(T)$ is well defined, then so is $f^*(S)$ for every 
positive closed current $S$ such that $S\leq T$.
\end{lemma}
\begin{proof}
We have $\lim\|T_n\|=\|T\|$. It follows that $\lim \|f^*(T_n)\|=d^p\|T\|
=\|f^*(T)\|$. On the other hand,
$f^*(T_n)\rightarrow f^*(T)$
in $\P^k\setminus f^{-1}(I^-)\cup I^+$ and $f^*(T)$ does not charge 
$f^{-1}(I^-)\cup I^+$. Hence $f^*(T_n)\rightarrow f^*(T)$ in $\P^k$.

We have $\|f^\star(S)\|\leq d^p\|S\|$,
$\|f^\star(T-S)\|\leq d^p\|T-S\|$
and $\|f^\star(T)\|=
d^p\|T\|$. It follows that $\|f^\star(S)\|= d^p\|S\|$.
Hence $f^*(S)$ is well defined.
\end{proof}

\begin{proposition} The operators
$f_*: \DSH^{k-p}(V)\rightarrow \DSH^{k-p}(V)$ and 
$f^*:\PC_p(V)\rightarrow\PC_p(V)$, $1\leq
p\leq s$, are well defined and are continuous.
We have
$(f^n)^*=(f^*)^n$,  $\|f^*(T)\|=d^p\|T\|$ and
$\langle f^*(T),\Phi\rangle = \langle T, f_*(\Phi)
\rangle$ for $T\in \PC_p(V)$ and $\Phi\in \DSH^{k-p}(V)$. 
\end{proposition}

\begin{proof} 
Let $\Phi\in\DSH^{k-p}(V)$. Using a partition of unity, we  can write 
$\Phi=\Phi^{(1)}+\Phi^{(2)}$ where $\Phi^{(1)}$
is a $\Linfty$ form with compact support in $V$ and $\Phi^{(2)}$ is a current with support in $\P^k\setminus I^+$.
By Definition 3.1, $f^{-1}:\P^k\setminus \overline U\rightarrow V$ and 
$f:\P^k\setminus I^+\rightarrow \P^k$ are holomorphic. Then $f_*(\Phi^{(1)})=(f^{-1})^*(\Phi^{(1)})$
and $f_*(\Phi^{(2)})$ are well defined. The first assertion follows, even $\Phi^{(1)}$ and 
$\Phi^{(2)}$ are not necessarily in $\DSH^{k-p}(V)$.

Consider now a smooth positive closed form $\Phi\in\DSH^{k-p}(V)$.
Recall that by Proposition 2.4, if $T$ is in $\PC_p(V)$, then $T$ and $f^\star(T)$ do not charge analytic sets.
We have 
$$\langle f^\star(T),\Phi\rangle = \langle T, f_*(\Phi)
\rangle_{\P^k\setminus I^-}:=\int_{\P^k\setminus I^-} T\wedge f_*(\Phi).$$
We next show that $\langle T, f_*(\Phi)
\rangle_{\P^k\setminus I^-} = \langle T, f_*(\Phi) \rangle$.
Let $W$ be a form, smooth outside $I^-$, such that
$\ddc W=f_*(\Phi)-m\omega^{k-p}$ and $\theta$ be a smooth function
supported in $U$, equal to 1 in a neighbourhood of $I^-$. Here, $m$ is the mass of $f_*(\Phi)$.
Define $\Psi:=\ddc(\theta W)+c\omega^{s-p}\wedge\Theta$, $c>0$ big enough.
Then $\Psi$ is positive closed,
$\supp(\Psi)\subset \P^k\setminus \overline V$ and $\Psi-f_*(\Phi)$
is smooth. This form $\Psi$ belongs to $\DSH^{k-p}(V)$.
We only need to show that
$\langle T, \Psi
\rangle_{\P^k\setminus I^-} = \langle T, \Psi \rangle$.
 
Let $u_n$ as in Lemma 2.6 but we replace $A$ by $I^-$.
We have
$$\langle T, \Psi
\rangle_{\P^k\setminus I^-} = \lim \langle T, u_n \Psi
\rangle = \langle T, \Psi\rangle$$
because $T\in\PC_p(V)$ and
$u_n\Psi\rightarrow \Psi$ in $\DSH^{k-p}(V)$.
So $\langle f^\star (T),\Phi\rangle = \langle T, f_*(\Phi) \rangle$
for $\Phi\in\DSH^{k-p}(V)$ smooth positive and closed. 

For $\Phi=\omega^{s-p}\wedge \Theta$, we get
$$\|f^\star(T)\|=\langle f^\star(T),\Phi\rangle
=  \langle T, f_*(\Phi)\rangle = d^p\|T\|.$$
The last equality follows from a regularization of the positive closed current 
$f_*(\Phi)$ and the properties:
$\|f_*(\Phi)\|=d^p$ and $T\in \PC_p(V)$. Hence $f^*(T)$
is well defined and equal to $f^\star(T)$. 

Assume now that $\Phi$ is a smooth positive form in $\DSH^{k-p}(V)$ 
non necessarily closed.
Using a regularization of $f_*(\Phi)$ in $\DSH^{k-p}(V)$, we get 
$$\langle f^* (T),\Phi\rangle=\langle T, f_*(\Phi)
\rangle_{\P^k\setminus I^-} \leq \langle T, f_*(\Phi) \rangle.$$
On the other hand, if $\Phi'\geq \Phi$ is a smooth closed form, we also have 
$$\langle f^*(T), \Phi'-\Phi\rangle=\langle T, f_*(\Phi'-\Phi)
\rangle_{\P^k\setminus I^-} \leq \langle T, f_*(\Phi'-\Phi) \rangle.$$
The equality $\langle f^*(T),\Phi'\rangle = \langle T, f_*(\Phi') \rangle$
implies that $\langle f^*(T),\Phi\rangle = \langle T, f_*(\Phi) \rangle$.
This also holds for $\Phi$ smooth non positive because we can write $\Phi$ as 
a difference of positive forms. 
From the first assertion of the Proposition, it follows that
the right hand side of the last equality 
is well defined
for every $\Phi\in\DSH^{k-p}(V)$ and depends
continuously on $\Phi$. This allows to extend $f^*(T)$ to a
continuous linear form on $\DSH^{k-p}(V)$. Hence $f^*(T)\in
\PC_p(V)$. The continuity of $f^*$ and the equality
$(f^n)^*=(f^*)^n$ are clear.
\end{proof}

\section{Convergence toward the Green currents}

Let $f$ be an $s$-regular 
birational map of algebraic degree $d\geq 2$ as in Definition 3.1.
Recall that the Green $(1,1)$-current $T_+:=\lim d^{-n}(f^n)^*(\omega)$ of $f$
has continuous local potentials 
in a neigbourhood of
$\P^k\setminus V$ \cite{Sibony2}.
Proposition 2.5 shows that $T_+^p$ is well defined
for $1\leq p\leq s$.
It belongs to $\PC_p(V)$.
Moreover,
we have $\lim d^{-np} f^{n*}(\omega^p) =
T_+^p$ in $\P^k\setminus\overline V$. The last property
follows from a uniform convergence of potentials
of $d^{-n} f^{n*}(\omega)$ (see 
\cite{Sibony2}). This is also reproved in Theorem 4.1.
We have $f^*(T_+^p)=d^p T_+^p$.

Let $I^+_n$ be the indeterminacy set of $f^n$. Define
 $U_\infty:=\cup_{n\geq 0} f^{-n}(U)\setminus I_n^+$.
In this section, we prove the following result 
which implies Theorem 1.2.

\begin{theorem} Let $f:\P^k\rightarrow\P^k$
be an $s$-regular birational map as above.
Then for every $p$, $1\leq p\leq s$, the following holds.
\begin{enumerate}
\item[{\rm 1.}] If $T\in \PC_p(V)$ is a positive closed current 
of mass $1$, then 
$d^{-pn}f^{n*}(T)$ converges in $U_\infty$ to $T^p_+$. Moreover,
every limit value of the sequence $d^{-np}f^{n*}(T)$
is in $\PC_p(V)$. The convergence is valid in the weak 
topology of $\PC_p(V)$. 
\item[{\rm 2.}] If $S$ is a positive closed $(p,p)$-current 
such that 
$S\leq T_+^p$ in $\P^k$, then $S=c T_+^p$ in $U_\infty$ where $c:=\|S\|$.
\end{enumerate}
\end{theorem}

\begin{proof}
1. Let $\Phi$ be a $(k-p,k-p)$-current in $\DSH^{k-p}(V)$. Write
$\ddc\Phi=\Omega=\Omega^+-\Omega^-$ where $\Omega^\pm$ are positive
closed $(k-p+1,k-p+1)$-currents supported in $\P^k\setminus
\overline V$.
Assume that $\|\Omega^\pm\|= 1$.
Define
$\Omega_n^\pm:=(f^n)_*(\Omega^\pm)$ and
$\Omega_n=\Omega_n^+-\Omega_n^-$ for $n\geq 0$.
They are supported in $U$ for $n\geq 1$ and we have 
$\|\Omega_n^\pm\|=d^{(p-1)n}$.

Let $\Phi_n^\pm$ be the solution of the equation
$\ddc\Phi_n^\pm=\Omega_n^\pm-d^{(p-1)n}
\omega^{s-p+1}\wedge \Theta$ given in
Proposition 2.1. The $\Phi^\pm_n$'s are negative $(k-p,k-p)$-forms, smooth on
$V$ and they satisfy
$\|\Phi_n^\pm\|_{\Linfty(V)}+
\|\Phi_n^\pm\|_\DSH\lesssim  d^{(p-1)n}$.
Define $\Phi_n:=\Phi_n^+-\Phi_n^-$, $\Psi_0:=\Phi-\Phi_0$ and
$\Psi_{n+1}:=f_*(\Phi_n)-\Phi_{n+1}$. The
forms $\Phi_n$ are smooth on $V$,
$\ddc\Phi_n=\Omega_n$ and $\|\Phi_n\|_\DSH\lesssim d^{(p-1)n}$ for $n\geq 1$.
By Proposition 3.5, $\|\Psi_n\|_\DSH\lesssim
d^{(p-1)n}$.
Since $\ddc\Psi_n=0$, we can associate to $\Psi_n$ a class
$b_n$ in $H^{k-p,k-p}(\P^k,\C)$.
We have $\|b_n\|\lesssim \|\Psi_n\|_{\Lone}\lesssim d^{(p-1)n}$.  

Since we assume
that $T\in\PC_p(V)$, Proposition 3.5 allows the following calculus 
\begin{eqnarray*}
\langle f^{n*}(T),\Phi\rangle & = &  \langle
f^{n*}(T),\Psi_0\rangle + \langle f^{n*}(T),\Phi_0\rangle\\ 
& = & \langle f^{n*}(T),\Psi_0\rangle 
+ \langle f^{(n-1)*}(T),f_*(\Phi_0)\rangle \\
& = & \langle f^{n*}(T),\Psi_0\rangle 
+ \langle f^{(n-1)*}(T), \Psi_1\rangle +
\langle f^{(n-1)*}(T),\Phi_1\rangle \\
& = & \langle f^{n*}(T),\Psi_0\rangle 
+ \langle f^{(n-1)*}(T), \Psi_1\rangle +
\langle f^{(n-2)*}(T),f_*(\Phi_1)\rangle.
\end{eqnarray*}
Using the equality $f_*(\Phi_n)=\Psi_{n+1}+\Phi_{n+1}$ we obtain
by induction that
\begin{eqnarray}
\langle f^{n*}(T),\Phi\rangle  & = &   
\langle f^{n*}(T),\Psi_0\rangle 
+ \langle f^{(n-1)*}(T), \Psi_1\rangle 
+ \cdots \nonumber\\
& & \cdots + \langle T, \Psi_n\rangle +
\langle T,\Phi_n\rangle.
\end{eqnarray}
Since $f^{n*}(T)$ is cohomologous
to $d^{pn}\omega^p$, using a regularization of $\ddc$-closed currents $\Psi_i$ in $\DSH^{k-p}(V)$, we get
$$\langle d^{-pn} f^{n*}(T),\Phi\rangle =
\int [\omega^p]\wedge (b_0+d^{-p}b_1+\cdots + d^{-pn}b_n) +
d^{-pn} \langle T,\Phi_n\rangle.$$  

Recall that $T\in \PC_p(V)$ and
$\|\Phi_n^\pm\|_{\Linfty(V)}+\|\Phi_n^\pm\|_\DSH
\lesssim d^{(p-1)n}$. It follows that $\lim d^{-pn}\langle T,\Phi_n\rangle =0$.  
The relations $\|b_n\|\lesssim d^{(p-1)n}$ imply that
\begin{eqnarray}
\lim \langle d^{-pn} f^{n*}(T),\Phi\rangle =
\int [\omega^p]\wedge c_\Phi & \mbox{ where } & 
c_\Phi:=\sum_{n\geq 0} d^{-pn}b_n.
\end{eqnarray}
Propositions 2.1 and 3.5 imply also that $c_\Phi$ depends continuously on
$\Phi\in \DSH^{k-p}(V)$. So, (3) implies that every limit value of
the sequence $d^{-pn} f^{n*}(T)$ belongs to $\PC_p(V)$.

Consider now a smooth real-valued $(k-p,k-p)$-form $\Phi$ 
supported in $U$. Observe that $\Phi$ can be written as a difference $\Phi_1-\Phi_2$ of 
negative forms supported in $U$ and that $\ddc\Phi_i+c(\omega^{s-p+1}\wedge \Theta)$ is
positive for $c>0$ big enough. 
It follows that $\Phi\in \DSH^{k-p}(V)$. By (3), 
$d^{-pn} f^{n*}(T)$ converges on $U$ to a current
which does not depend on $T$. Hence, $\lim d^{-pn} f^{n*}(T) =T_+^p$
on $U$ since this is true for  $T=\omega^p$ (and for $T=T_+^p$).
The relation $f^{n*}(T_+^p)=d^{np}T_+^p$ implies that
$\lim d^{-pn} f^{n*}(T) =T_+^p$
on $U_\infty$. 
\\

2. Let $c$ be the mass of $S$ and define $S_n:=d^{np}(f^n)_\star (S)$.
We have $S_n\leq T_+^p$. By Lemma 3.4, $f^{n*}(S_n)$ is well defined. 
From Proposition 2.4, $T_+^p$ has no mass on analytic sets. It follows that 
$f^{n*}(S_n)= d^{np}S$ since this holds out of an analytic set.
We also deduce that $\|S_n\|=c$.

Assume that $\Phi$ is smooth and supported in $U$.
Proposition 2.1 shows that $\Phi_j$ and $\Psi_j$ belong to the class ${\cal C}$ as in Lemma 2.6
for $A=\cup_{i\leq n} f^i(I^-)$.
Hence, we can apply Lemma 2.6 to
$R_i=T_+$ and to $(f^{n-j})^*(S_n)$. We get
$$\langle (f^{n-j+1})^*(S_n),\Psi_j\rangle = \langle (f^{n-j})^*,f_*(\Psi_j)\rangle$$
since these integrals can be computed out of the singularities of $f$, $\Psi_j$ and $f_*(\Psi_j)$. 
We can then apply 
(2) to $S_n-cT_+^p$.
Since $S_n-cT_+^p$ is cohomologous to 0, we get
\begin{eqnarray*}
d^{np}\langle S-cT_+^p,\Phi\rangle & = & 
\langle f^{n*}(S_n-cT_+^p),\Phi\rangle  =  
\langle S_n-cT_+^p,\Phi_n\rangle \\
& = &  
\langle S_n-cT_+^p,\Phi_n^+-\Phi_n^-\rangle.
\end{eqnarray*} 
The relations $S_n\leq T_+^p$
and $\Phi_n^\pm\leq 0$ imply 
that the last expression is dominated 
by a combination of $\langle T_+^p,\Phi_n^+\rangle$ and 
of $\langle T_+^p,\Phi_n^-\rangle$. Hence, since $T_+^p\in\PC_p(V)$ and 
$d^{-(p-1)n}\Phi_n^\pm$ belong to a compact set in $\DSH^{k-p}(V)$, we have
$$d^{np}|\langle S-cT_+^p,\Phi\rangle|\lesssim d^{(p-1)n}.$$
It follows that 
$\langle S-cT_+^p,\Phi\rangle=0$ for every smooth form $\Phi$
supported in $U$.
Hence, $S=cT_+^p$ on $U$.
In the same way, we show that $S_n=cT_+^p$ on $U$. The relations $f^{n*}(S_n)=d^{np}S$ and 
$f^{n*}(T_+^p)=d^{np} T_+^p$ imply that $S=cT_+^p$ on $f^{-n}(U)\setminus I_n^+$ 
for every  $n\geq 1$. 
\end{proof}

\begin{remark} \rm
The convergence in Theorem 3.1 is uniform on $T\in\PC_p(V)$
such that $|\langle T,\Phi\rangle|\leq
c(\|\Phi\|_{\Linfty(V)} + \|\Phi\|_\DSH)$, $c>0$, for every
$\Phi\in\DSH^{k-p}(V)$.  
\end{remark}

\small

Tien-Cuong Dinh and Nessim Sibony, 
Math\'ematique - B\^at. 425, UMR 8628, 
Universit\'e Paris-Sud, 91405 Orsay, France.\\ 
E-mails: (TienCuong.Dinh, Nessim.Sibony)@math.u-psud.fr

\end{document}